\documentclass[11pt]{article}
\usepackage{amsfonts,amsmath,amssymb,amsthm, amscd}
\usepackage[dvips]{epsfig,graphics}

\newtheorem{theorem}{Theorem}

\newtheorem{lemma}[theorem]{Lemma}

\theoremstyle{definition}

\newtheorem{example}[theorem]{Example}
\newtheorem{remark}[theorem]{Remark}

\def\<{{\langle}}
\def\>{{\rangle}}

\def\e{{\epsilon}}

\def\L{{\Lambda}}

\def\e{\epsilon}

\def\Z{\mathbb Z}
\def\C{\mathbb C}

\def\D{{\Delta}}

\begin{document}

\title{On Eigenvalues of Free Group Endomorphisms}

\author{Daniel S. Silver \and Susan G. Williams \\ {\em
{\small Department of Mathematics and Statistics, University of South Alabama}}}

\maketitle 

\begin{abstract} \noindent Let $\phi: F\to F$ be an endomorphism of a finitely generated free group, and let $H$ be a finite-index subgroup of $F$ that is invariant under $\phi$. The nonzero eigenvalues of $\phi$ are contained in the eigenvalues of $\phi$ restricted to $H$.      \end{abstract} 



\bigskip\bigskip

Let $F$ be a free group of finite rank, and let $\phi: F \to F$ be an endomorphism. By an 
{\it eigenvalue} of $\phi$ we mean an eigenvalue of the induced linear map
$\phi^{\rm ab} \otimes 1: F/[F,F]\otimes_\Z \C \to F/[F,F]\otimes_\Z \C.$ We denote the set of nonzero eigenvalues of $\phi$ by $\L(\phi)$. 

\begin{theorem}\label{main} Let $\phi: F \to F$ be an endomorphism of a finitely generated free group, and let $H$ be a finite-index subgroup of $F$ that is invariant under $\phi$. Then 
$\L(\phi) \subset \L(\phi\vert_H).$ \end{theorem}

\begin{remark} \label{rem} (1) Theorem \ref{main} answers a question attributed to A. Casson
 in \cite{bms} (F33): If $\phi$ is an automorphism of a finitely generated free group, does there exist a finite-index $\phi$-invariant subgroup $H$ such that every eigenvalue of $\phi\vert_H$ is a root of unity?\footnote{In \cite{best} the problem appears with ``every" replaced by ``some" (Question 12.16).  We thank Henry Wilton for bringing this to our attention.} 
A negative answer is provided by any automorphism $\phi$ possessing an eigenvalue that is not a root of unity. In particular, one can consider the automorphism $\phi$ of $F = \<a, b \mid \>$ mapping $a \mapsto b$ and $b \mapsto a b^2$. Its eigenvalues are $1 \pm \sqrt 2$. By Theorem \ref{main}, they are eigenvalues of the restriction $\phi\vert_H$, for any finite-index $\phi$-invariant subgroup $H \le F$. 
We return to this example at the conclusion of the paper. 

(2) By results of R. Bowen \cite{bowen}, the eigenvalues of $\phi\vert_H$ are bounded by the growth rate of $\phi$. If $\phi: F \to F $ is an endomorphism of a finitely generated group (not necessarily free), then its {\it growth rate} is defined to be 
$${\rm GR}(\phi) = \max_i \limsup_k ||\phi^k(x_i)||^{1/k},$$
where $\{x_i\}$ is a finite set of generators for $F$ and $||\cdot||$ is the word metric on $F$
(that is, $||y||$ is the minimum number of occurrences of $x_i^{\pm 1}$ in any word describing $y$). 
Bowen showed that ${\rm GR}(\phi)$ is finite and independent of the generator set. 
The well-known relevant facts are repeated here for the reader's convenience. 
\item{} $\bullet$ If $\phi_i: F_i \to F_i$, $i =1, 2$ are endomorphisms of finitely generated groups such that 
$\phi_1 h = h \phi_2$ for some epimorphism $h: F_1 \to F_2$, then ${\rm GR}(\phi_1) \ge {\rm GR}(\phi_2)$. In particular, ${\rm GR}(\phi) \ge {\rm GR}(\phi^{\rm ab})$ for any endomorphism $\phi$.

\item{}$\bullet$  The growth rate of an endomorphism of a finitely generated free abelian group is the maximum of the moduli of its eigenvalues.
\item{} $\bullet$ If $\phi$ is any endomorphism of a finitely generated group and $H$ is a $\phi$-invariant subgroup of finite index, then ${\rm GR}(\phi\vert_H) = {\rm GR}(\phi).$ 
\end{remark} 

We prove Theorem \ref{main}.

\begin{proof} We first consider the case that $\phi$ is injective. Let 
\begin{equation} \label{*} M(\phi) = \<F, x \mid x^{-1} g x = \phi(g),  \forall g \in F\> \end{equation}
and
\begin{equation}\label{**} M(\phi\vert_H) = \<H, x \mid x^{-1} g x = \phi(g),  \forall g \in H\>\end{equation}
be the algebraic mapping tori of $\phi$ and $\phi\vert_H$, respectively. Since $\phi$ is injective, 
each is presented as an ascending HNN decomposition (see \cite{ls}, for example). 
Elements of $M(\phi\vert_H)$ have the form  $h x^i$, where $h \in H$ and $ i \in \Z$.
By Britton's Lemma, such an element is trivial if and only if $h=1$ and $i = 0$. 
Consequently, the inclusion map $H \hookrightarrow F$ and assignment $x \mapsto x$
induces a natural embedding $M(\phi\vert_H) \hookrightarrow M(\phi)$. We regard $M(\phi\vert_H)$ as a subgroup of $M(\phi)$ via the embedding. 

The subgroup $M(\phi\vert_H)$ has finite index $m=[F:H]$ in $M(\phi)$. To see this, express $F$ as the union of disjoint cosets $a_1H, \ldots, a_m H$. Then 
$M(\phi) = \cup_{i \in \Z} F x^i = \cup_i \cup_j a_jHx^i \subset \cup_j a_j M(\phi\vert_H)$. If $a_j a_k^{-1}\in M(\phi\vert_H)$, then $a_j a_k^{-1}=hx^i$ for some $h\in H$, and Britton's Lemma gives $a_j a_k^{-1}\in H$.

Let $X_{M(\phi)}$ be a standard 2-dimensional cell complex with $\pi_1 X \cong M(\phi)$, corresponding to the presentation (\ref{*}). Consider the commutative diagram of covering space maps:

$$\begin{matrix} X_H & {\buildrel q \over \longrightarrow} & X_F \\

\Big\downarrow &  & \Big\downarrow\\

X_{M(\phi\vert_H)} &  {\buildrel p \over \longrightarrow} & X_{M(\phi)} 

\end{matrix}$$

\noindent Here $X_{M(\phi\vert_H)}$ is the connected cover of degree $[F:H]$ with ${\rm im}\ p_* = M(\phi\vert_H)$.  Also,
$X_F$ and $X_H$ are regular covers with fundamental groups $F$ and $H$, respectively, corresponding to the homomorphism $\e: M(\phi) \to \Z = \<t \mid \>$ and its restriction to 
$M(\phi\vert_H)$ mapping $F$ trivially and sending $x$ to $t$. Then $p$ lifts to a covering map $q$ with ${\rm im}\ q_* = H$.

The homology group $H_1(X_F ; \C)$ is a finitely generated module over $\C[t, t^{-1}]$. Its module order is $\D(t) = \det(t I - \phi^{\rm ab} \otimes 1)$. Similarly, $H_1(X_H; \C)$ has module order 
$\tilde \D(t) = \det(t I - (\phi\vert_H)^{\rm ab} \otimes 1)$. Both $\D$ and $\tilde \D$ can be readily computed from the presentations (\ref{*}) and (\ref{**}) using Fox calculus (see \cite{tur}, for example). 
Each is well defined up to the units of $\C[t, t^{-1}]$, monomials with nonzero complex coefficients.

The covers $X_F$ and $X_H$ acquire cell complex structures from $X$. The transfer homomorphism 
$\tau: H_1(X_F; \C) \to H_1(X_H; \C)$ is a module monomorphism (see \cite{hatcher}). Since $H_1(X_F; \C)$ embeds as a submodule of $H_1(X_H; \C)$, the order $\D$ divides $\tilde \D$ (see \cite{milnor}). Since the nonzero roots of $\D$ (resp. $\tilde \D$) are the eigenvalues of $\phi$ (resp. $\phi\vert_H$), we have $\L(\phi) \subset \L(\phi\vert_H)$. 

For general endomorphisms, we adapt an approach of I. Kapovich  \cite{kap}. Consider the nested sequence
$${\rm ker}\ \phi \subset {\rm ker}\ \phi^2 \subset \cdots \subset {\rm ker}\ \phi^i \subset \cdots.$$
As in \cite{kap}, \cite{sela}, there exists an integer $k >0$ such that ${\rm ker}\ \phi^k = {\rm ker}\ \phi^n$ whenever $n \ge k$. Let $K = {\rm ker}\ \phi^k$, called here the {\it eventual kernel} of $\phi$.  There is an induced injection $\bar \phi: F/K \to F/K$. The group $F/K$ is isomorphic to $\phi^k(F)$, and hence free and of finite rank. It is shown in \cite{kap} that $\phi$ and $\bar\phi$ have isomorphic mapping tori, a fact that we do not use here.  We will want the following lemma, proved at the end.

\begin{lemma}\label{lem} $\L(\bar \phi) = \L(\phi)$. \end{lemma}

Since ${\rm ker}\ (\phi\vert_H)^i = {\rm ker}\ \phi^i \cap H$ for each $i > 0$, the eventual kernel of $\phi\vert_H$ is $H \cap K$, and $H/H\cap K$ is isomorphic to $\bar H=\phi^k(H)$.  Writing $F = a_1H \cup \cdots \cup a_mH$ as above, we see $\phi^k(F)=\phi^k(a_1)\bar H\cup\cdots\cup \phi^k(a_m)\bar H$, so $[\phi^k(F):\phi^k(H)]$ is finite.

%

The proof for injective $\phi$ together with Lemma \ref{lem} yield
$$\L(\phi) = \L(\bar \phi) \subset  \L(\bar \phi\vert_H) =  \L(\phi\vert_H).$$  \end{proof}

We prove Lemma \ref{lem}. 

\begin{proof}  The natural short exact sequence 
$$1 \to K \to F \to F/K \to 1$$
splits since $F/K$ is free. Then $F$ is a semidirect product 
$K \rtimes_\theta F/K$, where $\theta: F/K \to {\rm Aut}\ K$, $y \mapsto \theta_y$, is a free group action. More explicitly, $F$ is isomorphic to the free product of $K$ and $F/K$ modulo the relations 
$y u y^{-1} = \theta_y(u)$, where $y\in F/K$ and $u \in K$. The abelianization
$F^{\rm ab}$ is isomorphic to $K^{\rm ab}/M \times (F/K)^{\rm ab} $, where $M$ is the subgroup generated by the elements $[u] - \theta^{\rm ab}_y([u])$. Here $\theta^{\rm ab}_y$ is the automorphism of $K^{\rm ab}$ induced by $\theta_y$, and $[u]$ is an arbitrary element of $K^{\rm ab}$. 

Since $K$ is $\phi$-invariant, $K^{\rm ab}$ is $\phi^{\rm ab}$-invariant. Consider the $\phi^{\rm ab}$-image of
$[u] - \theta^{\rm ab}_y([u])$. It can be written as $[\phi(u)]-[\phi(\theta_y(u))]$, or equivalently $[\phi(u)]-[\phi(y) \phi(u) \phi(y)^{-1}]$. Since $\phi(y)$ is represented by $v y'$ for some 
$v \in K, y' \in F/K$, we can write this element as $[\phi(u)]-[v y' \phi(u) y'^{-1} v^{-1}]$. But $[v y' \phi(u) y'^{-1} v^{-1}]= [y' \phi(u) y'^{-1}] = [\theta_{y'} (\phi(u))]$. Hence the $\phi^{\rm ab}$-image of 
$[u] - \theta^{\rm ab}_y([u])$ is $[\phi(u)]- [\theta_{y'}(\phi(u))]$, which is in $M$. 

We have shown that $M$ is invariant under $\phi^{\rm ab}$. Hence  $\phi^{\rm ab}$ induces 
a homomorphism of $K^{\rm ab}/M$.
The linear map $\phi^{\rm ab}: K^{\rm ab} \to K^{\rm ab}$ has a matrix representation of the form $$\begin{pmatrix} A & 0 \\ * & B \end{pmatrix},$$ where  $A$ represents $\bar \phi^{\rm ab}$ and $B$ represents $\phi^{\rm ab}$ restricted to $K^{\rm ab}/M$. Since $B$ is clearly nilpotent, the sets of nonzero eigenvalues of $\bar \phi$ and $\phi$ are equal. 

\end{proof}

\begin{example} We return to the example in Remark \ref{rem}. The linear map $\phi^{\rm ab}$ has a matrix representation 
$$\begin{pmatrix} 0 & 1 \\ 1 & 2 \end{pmatrix}.$$ Here $\L(\phi) = \{ 1+\sqrt 2, 1- \sqrt 2\}$. 
The subgroup $H < F$ generated by $a^2, b^2, ab$ has index equal to 2 and is $\phi$-invariant. 
Since $\phi(a^2) = b^2,\ \phi(b^2) = (ab)b^2 (ab)^{-1} a^2 b^2,\ \phi(ab) = b^2 (ab)^{-1}a^2 b^2$, the map $(\phi\vert_H)^{\rm ab}$ has a matrix representation 
$$\begin{pmatrix} 0&1&1\\1&2&2\\0 &0& -1\end{pmatrix}.$$
Hence $\L(\bar \phi) = \{1+\sqrt 2, 1- \sqrt 2, -1\}$. 
\end{example}



\begin{thebibliography}{10}
 
 \bibitem{bms} G. Baumslag, A.G. Myasnikov and V. Shpilrain, {\it Open problems in combinatorial group theory. Second edition.}, in Cont. Math. {\bf 296}, American Mathematical Society, Providence, RI, 2002. Online version: www.sci.ccny.cuny.edu/~shpil/gworld/problems/oproblems.html
 
 \bibitem{best} M. Bestvina, {\it Questions in Geometric Group Theory} (updated 2004), www.math.utah.edu/$\sim$bestvina/eprints/questions-updated.pdf.
 
 \bibitem{bowen} R. Bowen, {\it Entropy and the fundamental group}, in The structure of attractors in dynamical systems (Proc. Conf., North Dakota State Univ., Fargo, ND, 1977), 21--29, Lecture Notes in Math., {\bf 668}, Springer, Berlin, 1978.
 
 \bibitem{hatcher} A. Hatcher, Algebraic Topology, Cambridge University Press, Cambridge, 2002.
 
 \bibitem{kap} I. Kapovich, {\it A remark on mapping tori of free group endomorphisms}, arXiv:math/0208189v1
 
 \bibitem{ls} R.C. Lyndon and P.E. Schupp, Combinatorial group theory, Springer-Verlag, Berlin, 2001.
 
 \bibitem{milnor} J.W. Milnor,  {\it Infinite cyclic covers}, in Conf. Topology of Manifolds 1968 (ed. J.G. Hocking), 115 -- 133, Prindle, Weber and Schmidt, Boston, 1968. 
 
 \bibitem{sela} Z. Sela, {\it Endomorphisms of hyperbolic groups. I. The Hopf property}, {\sl Topology\ \bf 38}(1999), 301--321. 
 
 \bibitem{tur} V. Turaev, Introduction to Combinatorial Torsions, Birkh\"auser Verlag, Basel, 2001. 

 
 
 
 \end{thebibliography}
\end{document}